\documentclass[letterpaper, 10 pt, conference]{ieeeconf}
\usepackage{blkarray}
\usepackage{dsfont}
\usepackage{bm}
\usepackage{comment}
\usepackage[cmex10]{amsmath}
\usepackage{amssymb}
\usepackage{gensymb}
\usepackage{epsfig,floatflt}
\usepackage{psfrag,graphicx}
\usepackage[colorlinks]{hyperref}
\usepackage{cite}
\usepackage{epstopdf}
\usepackage{subcaption}
\usepackage[labelformat=parens,labelsep=quad,skip=3pt]{caption}
\IEEEoverridecommandlockouts

\newcommand{\bx}{\bm{\bar{x}}}
\newcommand{\bu}{\bm{\bar{u}}}

\newcommand{\IR}{\mathbb{R}}

\newtheorem{thm}{\textbf{Theorem}}
\newtheorem{prob}{\textbf{Problem}}

\newtheorem{lem}{\textbf{Lemma}}
 
\newtheorem{defn}{\textbf{Definition}}

\newtheorem{asm}{\textbf{Assumption}}
\newtheorem{rem}{\textbf{Remark}}

\DeclareSymbolFont{extraup}{U}{zavm}{m}{n}
\DeclareMathSymbol{\vardiamond}{\mathalpha}{extraup}{87}
\newcommand*{\qedprob}{\hfill\ensuremath{\triangledown}}%

\begin{document}
	\title{Optimal Control using Composite Bernstein Approximants}
	\author{Gage MacLin$^{1}$, Venanzio Cichella$^{1}$, Andrew Patterson$^{2}$, Michael Acheson$^{2}$, and Irene Gregory$^{2}$.
		\thanks{This work was supported by NASA and ONR.}
		\thanks{$^{1}$ Gage MacLin and Venanzio Cichella are with the Department of Mechanical Engineering, University of Iowa, Iowa City, IA 52240 {\tt\small \{gage-maclin, venanzio-cichella\} @uiowa.edu}}
        \thanks{$^{2}$ Andrew Patterson, Michael Acheson, and Irene Gregory are with the NASA Langley Research Center, Hampton, Virginia, VA 23681 {\tt\small \{andrew.patterson, michael.j.acheson, irene.m.gregory\} @nasa.gov}}
	}
	\maketitle
	\begin{abstract}
		In this work, we present composite Bernstein polynomials as a direct collocation method for approximating optimal control problems. An analysis of the convergence properties of composite Bernstein polynomials is provided, and beneficial properties of composite Bernstein polynomials for the solution of optimal control problems are discussed. The efficacy of the proposed approximation method is demonstrated through a bang-bang example. Lastly, we apply this method to a motion planning problem, offering a practical solution that emphasizes the ability of this method to solve complex optimal control problems.
	\end{abstract}
	

	\IEEEpeerreviewmaketitle


	\section{Introduction}
\label{sec:Introduction}
Optimal control problems are often too complex to be solved analytically, thus requiring the application of numerical methods to find a solution. Numerical methods for solving optimal control problems are typically separated into two primary classifications: direct and indirect methods. Indirect methods involve the transformation of the optimal control problem into a boundary-value problem, resulting in a system of differential equations. These methods rely on the calculus of variations and Pontryagin's maximum principle, often resulting in solutions that can be harder to implement for complex problems. Alternatively, direct methods discretize the continuous optimal control problem into a nonlinear programming (NLP) problem, allowing the solution to be found via collocation and pseudospectral (PS) methods.

A common approach to solving optimal control problems using direct methods is to employ PS methods, which have been used to solve a wide range of optimization problems 
\cite{karpenko2012first,ross2006pseudospectral}, due in large part to their spectral (exponential) rate of convergence. Various collocation methods have been introduced, including Legendre-Gauss-Lobatto (LGL) PS for finite time-horizon problems \cite{elnagar1995pseudospectral, ross2004legendre}, as well as Legendre-Gauss (LG) PS \cite{benson2005gauss, huntington2007advancement}
and Legendre-Gauss-Radau (LGR) PS \cite{kameswaran2008convergence, garg2011direct} for infinite time-horizon problems. Although PS methods have gained extensive adoption solving real-world optimal control problems, they are not without their own set of limitations.

The primary drawbacks of PS methods lie in the discretization and approximation of the optimal control problem \cite{ross2012review}. When discretizing the state and input of the optimal control problem, the constraints are enforced only at the specific discretization points, with no guarantees on constraint satisfaction between each point. Additionally, PS methods struggle to approximate discontinuous solutions due to the Gibbs phenomenon, which can be observed as oscillations around a jump discontinuity. There are methods of reducing the effects of the Gibbs phenomenon, though the most powerful of which require that the location of the discontinuities are known \emph{a priori} \cite{sarra2009edge}, a condition that is not feasible for many optimal control problems. Confronted with these challenges, researchers have been motivated to explore more robust alternatives that are not plagued by these limitations.

One alternative approach is Bernstein polynomial approximation, which has many properties that aid in the accurate solution of optimal control problems \cite{kielas2022bernstein}. Due to geometric properties of Bernstein polynomials, constraint satisfaction can be guaranteed throughout the solution, not solely at the discretization points. Furthermore, Bernstein polynomials do not suffer from the Gibbs phenomenon, allowing for smooth approximations of discontinuous solutions. With these properties, Bernstein polynomials are able to approximate the solution for the whole time horizon, without inducing oscillations at jump discontinuities. Composite Bernstein polynomials, a series of Bernstein polynomials segmented together at various knotting locations along the time horizon, enjoy these same properties as Bernstein polynomials. However, unlike single Bernstein polynomials, composite Bernstein polynomials are able to accurately approximate both smooth and discontinuous solutions and converge to the solution at a faster rate than single Bernstein polynomials, as presented in this work. Due to the introduction of these knotting locations, composite Bernstein polynomials are capable of accurately approximating discontinuous solutions, such as bang-bang control inputs, which are a common occurrence among optimal control problems. This paper extends our previous work on Bernstein polynomial-based optimal control to composite Bernstein polynomials, and proposes a new direct method to solve optimal control problems using these approximants.

This paper is structured as follows: in Section \ref{sec:ProblemStatement} we present the Bolza-type optimal control problem. Section \ref{sec:CBA} introduces composite Bernstein polynomials and their convergence properties. Section \ref{sec:ProposedDirectMethod} presents the proposed direct method with the discretization of the Bolza-type optimal control problem. In Section \ref{sec:KnottingMethod}, a method to determine the quantity of knots is discussed. Numerical examples are provided in Section \ref{sec:NumericalResults}, including a collision avoidance motion planning problem. Finally, conclusions are presented in Section \ref{sec:Conclusion}.

\section{Problem Statement}\label{sec:ProblemStatement}
Consider the Bolza-type optimal control problem:

\begin{prob}[Problem $P$]
\label{prob:continuous}
	\begin{equation} \label{eq:costfunc}
	\begin{split}  
	& \min_{\bm{x}(t),\bm{u}(t),{{t}_f}}  I(\bm{x}(t),\bm{u}(t),{{t}_f}) = \\  &  E(\bm{x}(0),\bm{x} 
 (t_f),t_f)+ \int_0^{t_f} F(\bm{x}(t),\bm{u}(t))dt \, 
	\end{split} 
	\end{equation}
subject to 
	\begin{align}
	& \dot{\bm{x}}(t) = \bm{f}(\bm{x}(t),\bm{u}(t))\, , \quad \forall t \in[0,t_f], \label{eq:dynamicconstraint} \\
	& \bm{e}(\bm{x}(0),\bm{x}(t_f),{{t}_f}) = \bm{0} \, , \label{eq:equalityconstraint} \\
	& \bm{h}(\bm{x}(t),\bm{u}(t)) \leq \bm{0} \, , \quad \forall t\in [0,t_f] \label{eq:inequalityconstraint} \, ,
	\end{align}
	where $I:\IR^{n_x} \times \IR^{n_u} \times \IR \to \IR$, $E:\IR^{n_x}\times \IR^{n_x} \times \IR \to \IR$, $F:\IR^{n_x} \times \IR^{n_u} \to \IR$, $\bm{f}:\IR^{n_x}\times \IR^{n_u}\to \IR^{n_x}$, $\bm{e}:\IR^{n_x} \times \IR^{n_x} \times \IR \to \IR^{n_e}$, and $\bm{h}:\IR^{n_x} \times \IR^{n_u}  \to \IR^{n_h}$.
	\qedprob
\end{prob}

For Problem $P$, $I$ represents the Bolza-type cost function, $E$ is the initial point and end point cost, and $F$ is the running cost; Equation \eqref{eq:dynamicconstraint} governs the dynamics of the system, Equation \eqref{eq:equalityconstraint} imposes boundary conditions, and Equation \eqref{eq:inequalityconstraint} ensures compliance with the inequality constraints of the system.

\section{Composite Bernstein Approximation}\label{sec:CBA}
In this section, we describe the process of function approximation through composite Bernstein polynomials, and present key convergence results.
The Bernstein basis for polynomials of degree $N$ over the domain $I_k = [t_{k-1},t_{k}]$ is defined as
\begin{equation*}
    b_{j,N}^{[k]}(t)= \binom{N}{j} \frac{(t-t_{k-1})^j(t_{k}-t)^{N-j}}{(t_k-t_{k-1})^N}
\end{equation*}
for $j=0,...,N$, where $\binom{N}{j}$ is the combination of $N$ and $j$.
An $N$th order Bernstein polynomial ${x}_N^{[k]}(t)$ defined over the domain $I_k$ is a linear combination of $N+1$ Bernstein basis polynomials of order $N$, i.e.,
\begin{equation} \label{eq:singlebpoly}
    {x}_N^{[k]}(t) = \sum_{j=0}^N  \bar{x}_{j,N}^{[k]} b_{j,N}^{[k]}(t)  \, , \qquad \forall t \in I_k, \forall k={1,...,K}.
\end{equation}
where $\bar{{x}}_{j,N}^{[k]}$ for $j=0,...,N$ are referred to as Bernstein coefficients (or control points). 
Given time \textit{knots} \(t_k\), \(k = 0, \ldots , K\), with \(t_{0} < t_{1}<\ldots<t_K\),
a \emph{composite Bernstein polynomial} is defined as follows:
\begin{equation} \label{eq:cbp}
     {x}_M(t)= x_N^{[k]}(t), \quad \forall t\in I_k, \forall k={1,...,K}.
\end{equation}

The derivative of a composite Bernstein polynomial ${x}_M(t)$ is defined for the closed set $[0,t_f]$, and is computed as
\begin{equation}
        \dot{{x}}_M(t) = \sum\limits_{j=0}^{N-1} \sum\limits_{i=0}^N \bar{{x}}_{i,N}^{[k]}D_{i,j}^{[k]}b_{j,N-1}^{[k]}(t), \quad  t \in I_k
\end{equation}
where $D_{i,j}^{[k]}$ is the $(i,j)$th entry of the Bernstein differentiation matrix
\begin{equation*}
\begin{split}
\bm{D}_{N-1}^{[k]} =
\begin{bmatrix}
- \frac{N}{t_k-t_{k-1}}  &  0 & \ldots & 0 \\
\frac{N}{t_k-t_{k-1}} &  - \frac{N}{t_k-t_{k-1}} & \ddots & \vdots \\ 
0 &  \ddots & \ddots & 0 \\
\vdots & \ddots & \frac{N}{t_k-t_{k-1}}  & -\frac{N}{t_k-t_{k-1}} \\
0 &  \ldots & 0 & \frac{N}{t_k-t_{k-1}}
\end{bmatrix} \\
\in\mathbb{R}^{(N+1)\times N}, \; \forall k=1,...,K.
\end{split}
\end{equation*}
The definite integral of a composite Bernstein polynomial $x_M(t)$ is computed as
\begin{equation}
    \int_0^{t_f}{x}_M(t)dt= 
    \sum_{k=1}^K w^{[k]} \sum\limits_{j=0}^N \bar{{x}}_{j,N}^{[k]}
\end{equation}
where $w^{[k]}=\frac{t_k-t_{k-1}}{N+1} , \; \forall k = 1, \ldots, K$.

Composite Bernstein polynomials can be used to approximate functions.
\begin{defn} \label{def:approximant}
Let \(x(t)\) be a function defined over \([0,t_f]\). Let the time knots satisfy \(t_0 = 0\), \(t_K = t_f\), and $|t_k-t_{k-1}| \le \frac{C_t}{K}$ for all \(k=1,\ldots,K\), where $C_t>0$ is independent of $K$. Let
\begin{equation}\label{eq:tkj}
t_{k,j} = t_{k-1} + \frac{j}{N}(t_{k} - t_{k-1}).
\end{equation}
Then, the \emph{composite Bernstein approximant} for \(x(t)\) is a composite Bernstein polynomial \({x}_M(t)\) with coefficients \(\bar{x}_{j,N}^{[k]} = x(t_{k,j})\).
\end{defn}

Next, we detail convergence results for composite Bernstein polynomials.

\begin{lem} Let \(x(t) \in \mathcal{C}^2([0,t_f])\). Let $x_M(t)$ be the composite Bernstein polynomial approximation of \(x(t)\). The following bound holds:   
\begin{equation*}
\begin{split}
\Vert{} {x}_M(t) - {x}(t) \Vert  \leq \frac{A}{K^2N}
\end{split}
\end{equation*}
for all $t \in [0,t_f]$, where
\begin{equation*}
    A = \frac{C_t^2}{8}\max_{\tau\in[t_{k-1},t_k],k=1,...,K}|\ddot{x}(\tau)|.
\end{equation*}
\label{lem:function}
\end{lem}
\textbf{Proof:} The outline of the proof of Lemma \ref{lem:function} is given in Appendix \ref{app.lem:function}.

\begin{lem}\label{lem:derivative} Let $x(t)\in \mathcal{C}^{r+2}(0,t_f)$ for some r $\in \mathbb{Z}^+$. Let $x_M(t)$ be the composite Bernstein polynomial approximation of \(x(t)\). The following bound holds:
\begin{equation*}
\Vert x_M^{(r)}(t)-x^{(r)}(t)\Vert\le\frac{B_1}{K^2N}+\frac{B_2}{KN}+\frac{B_3}{N}
\end{equation*}
where $x^{(r)}(t)$ denotes the rth derivative of $x(t)$, and $B_1$, $B_2$, and $B_3$ are constants independent of both $K$ and $N$:
\begin{equation}\label{eq:lemma2results}
\begin{split}
    B_1=\frac{C_t^2(r^2+1)}{2}\max_{\tau\in[t_{k-1},t_k],k=1,...,K}\Vert x^{(r+2)}(\tau)\Vert, \\
    B_2=\frac{C_tr^2}{2}\max_{\tau\in[t_{k-1},t_k],k=1,...,K}\Vert x^{(r+1)}(t_{k-1})\Vert, \\
    B_3=\frac{r(r-1)}{2}\max_{\tau\in[t_{k-1},t_k],k=1,...,K}\Vert x^{(r)}(t_{k-1})\Vert.
\end{split}
\end{equation}
\end{lem}
\textbf{Proof:} The outline of the proof of Lemma \ref{lem:derivative} is given in Appendix \ref{app.lem:derivative}.

\begin{lem}\label{lem:integral}
Let $x_M(t)$ be the composite Bernstein polynomial approximation of $x(t) \in \mathcal{C}^2$. The following bound holds:
\begin{equation*}
\begin{split}
    \left\Vert\int^{t_f}_{0}x_M(t)dt-\int^{t_f}_{0}x(t)dt\right\Vert
    \le\frac{C}{K^2N}
\end{split}
\end{equation*}
where
\begin{equation*}
    C = \frac{C_t^3}{8}\max_{\tau\in[t_{k-1},t_k]}|\ddot{x}(\tau)|.
\end{equation*}
\end{lem}
\textbf{Proof:} The proof of Lemma \ref{lem:integral} is given in Appendix \ref{app.lem:integral}.

Notice that the convergence rate of composite Bernstein polynomials is quadratic w.r.t. $K$, in contrast to the linear convergence rate observed with single Bernstein polynomials (w.r.t. $K$, where \(K=1\)), see also \cite{cichella2018bernstein}. Building upon this enhanced convergence rate, this paper extends our previous work on Bernstein polynomial-based optimal control to composite Bernstein polynomials. 

\section{Proposed Direct Method} \label{sec:ProposedDirectMethod}
Here we formulate a discretized version of Problem $P$, referred to as Problem $P_{M}$. Where $M+1=K(N+1)$ is the number of control points of the composite Bernstein polynomial,  with $K$ denoting the number of polynomials in the approximant and $N$ denoting the order of each polynomial. 

We approximate the states and control inputs of Problem $P$ with composite Bernstein polynomials, with each individual polynomial defined as:
\begin{equation}
\begin{split}
    \bm{x}_N^{[k]}(t)=\sum_{j=0}^N \bar{\bm{x}}_{j,N}^{[k]}b_{j,N}(t), \quad t \in [t_{k-1},t_k] \\
    \bm{u}_N^{[k]}(t)=\sum_{j=0}^N \bar{\bm{u}}_{j,N}^{[k]}b_{j,N}(t), \quad t \in [t_{k-1},t_k] 
\end{split}
\end{equation}
\(\forall k = 1,...,K\). 

The composite Bernstein polynomials $\bm{x}_M(t)$ and $\bm{u}_M(t)$, defined in the same form as Equation \eqref{eq:cbp}, approximate $\bm{x}(t)$ and $\bm{u}(t)$
with $\bm{x}_M : [t_0,t_K] \to \mathbb{R}^{n_x}$ and $\bm{u}_M : [t_0,t_K] \to \mathbb{R}^{n_u}$, where $t_0=0$ and $t_K$ is an approximation of the optimal final mission time $t_f$. Let 
\[\bar{\bm{x}}_M = [\bar{\bm{x}}_{0,N}^{[1]}, \ldots, \bar{\bm{x}}_{N,N}^{[1]}, \ldots,\bar{\bm{x}}_{0,N}^{[K]},\ldots,\bar{\bm{x}}_{N,N}^{[K]}],  \]
\[\bar{\bm{u}}_M = [\bar{\bm{u}}_{0,N}^{[1]}, \ldots, \bar{\bm{u}}_{N,N}^{[1]},\ldots,\bar{\bm{u}}_{0,N}^{[K]},\ldots,\bar{\bm{u}}_{N,N}^{[K]}],  \]
i.e., \(\bar{\bm{x}}_M\in \mathbb{R}^{n_x\times (M+1)}\) and \(\bar{\bm{u}}_M \in \mathbb{R}^{n_u\times (M+1)}\) are vector-valued functions that comprise all the control points from each Bernstein polynomial. 
Let $\bar{\bm{t}}_K$ be a vector of the time knots between each polynomial, i.e.,  $\bar{\bm{t}}_K=[t_0,t_1,...,t_K] \in \mathbb{R}^{K+1}$.
Then, Problem $P_{M}$ can be defined as follows:
\begin{prob}[Problem $P_{M}$] 
\label{prob:discrete}
\begin{multline}    \min_{\bar{\bm{x}}_M,\bar{\bm{u}}_M,t_K} I_M(\bar{\bm{x}}_M,\bar{\bm{u}}_M,\bar{\bm{t}}_K) = E(\bar{\bm{x}}_{0,M},\bar{\bm{x}}_{M,M},t_K) \\ +\sum_{k=1}^Kw^{[k]}\sum_{j=0}^{N}F(\bar{\bm{x}}_{j,N}^{[k]},\bar{\bm{u}}_{j,N}^{[k]})
\end{multline}
subject to
\begin{multline} \label{eq:discdynconstraint}
    \left\Vert\sum_{i=0}^{M} \bar{\bm{x}}_{i,M}\bm{D}_{i,j}^M-\bm{f}(\bar{\bm{x}}_{j,M},\bar{\bm{u}}_{j,M})\right\Vert\le \delta_P^M, \\ \forall j=0,...,M
\end{multline}
\begin{equation} \label{eq:discequalityconstraint}
    \bm{e}(\bar{\bm{x}}_{0,M},\bar{\bm{x}}_{M,M},t_K)=0,
\end{equation}
\begin{equation} \label{eq:discinequalityconstraint}
\bm{h}(\bar{\bm{x}}_{j,M},\bar{\bm{u}}_{j,M}) \le 0, \; \forall j = 0,...,M
\end{equation}
\begin{equation} \label{eq:disccontinuityx}
    \bar{\bm{x}}_{N,N}^{[k]}-\bar{\bm{x}}_{0,N}^{[k+1]} = 0, \; \forall k=1,...,K-1
\end{equation}
\begin{equation} \label{eq:disccontinuityt}
    t_K > t_{K-1} > \ldots > t_0 = 0, 
\end{equation}
where $w^{[k]} = \frac{t_k-t_{k-1}}{M+1}, \; \forall k =1,...,K$,  $\bar{\bm{x}}_{j,M}$ is the $j$th element of $\bar{\bm{x}}_M$, $\bar{\bm{u}}_{j,M}$ is the $j$th element of $\bar{\bm{u}}_M$, and $\delta_P^M$ is a relaxation bound equal to a small positive number that depends on \(M\) and converges uniformly to \(0\), i.e.,
\(
    \lim_{M\rightarrow{\infty}}\delta_P^M=0. 
\)
Finally, $\bm{D}_{i,j}^M$ is the $(i,j)$th entry of the differentiation matrix $\bm{D}_M\in \mathbb{R}^{(M+1) \times (M+1)}$,
\[
\bm{D}_M =
\text{blkdiag}\left({\bm{D}_{N-1}^{[1]}\bm{E}_{N-1}^N, \ldots, \bm{D}_{N-1}^{[K]}\bm{E}_{N-1}^N} \right)
 \]
where \(\text{blkdiag}()\)  is defined as the block diagonal operator, and $\bm{E}_{N-1}^N$ is the degree elevation matrix of order $N-1$ to $N$, see \cite{kielas2022bernstein}. Differentiation of Bernstein polynomials results in a polynomial of \(N-1\) order, degree elevation returns $D_{N-1}$ to a square matrix.
\end{prob}


The following results can be stated regarding the feasibility and consistency of the proposed method. 

\begin{thm} \label{thm:feasibility}
Assume that Problem \ref{prob:continuous} is feasible, and the solution satisfies \(\bm{x} \in \mathcal{C}^3\), \(\bm{u} \in \mathcal{C}^2\). Let the functions \(\bm{f}\) and \(\bm{h}\) in Problem \ref{prob:continuous} be Lipschitz with respect to their arguments. Then, there exist order of approximations \(K^*\) and \(N^*\), and for any approximation orders \(K \geq K^*\) and \(N \geq N^*\) there exist Bernstein coefficients \(\bar{\bm{x}}_M\), \(\bar{\bm{u}}_M\) and final time \(t_K\) that constitute a feasible solution to Problem \ref{prob:discrete}.
\end{thm}

\textbf{Proof:}
Let $\bm{x}(t)$ and $\bm{u}(t)$ be a solution for Problem $P$, which exists by assumption.
Define 
\begin{equation} \label{eq:coeffdef}
\bx_{j,N}^{[k]} = \bm{x}(t_{k,j}) , \, \bu_{j,N}^{[k]} = \bm{u}(t_{k,j})  , \, 
\end{equation}
for all \( k \in \{1,\ldots,K\} \,\), \( j \in \{0,\ldots,N\} \,\), where $t_{k,j}$ is defined in Equation \eqref{eq:tkj}.
Under the assumption \(\bm{x}\in \mathcal{C}^3, \, \bm{u} \in \mathcal{C}^2\), Lemma \ref{lem:function} implies that
\begin{equation} \label{eq:unifconvxnprime}
\begin{split}
& ||\bm{x}_M(t) - \bm{x}(t)|| \leq A_x/K^2N  \, ,  \\
& ||\bm{u}_M(t) - \bm{u}(t)|| \leq A_u/K^2N  \, , \\
\end{split}
\end{equation}
for all $t\in [0,t_f]$, where $\bm{x}_M(t)$ and $\bm{u}_M(t)$ are computed as in Equation \eqref{eq:cbp} with coefficients given by Equation \eqref{eq:coeffdef}, and $A_x,A_u>0$ are independent of \(K\) and \(N\).
Next, we show that the above polynomials satisfy the constraints in Problem \ref{prob:discrete}. For all \(k=1,\ldots,K\) and \(j=0,\ldots,N\), the left hand side of Equation \eqref{eq:discdynconstraint} gives 
\begin{equation*}
\begin{split}
& ||\sum_{i=0}^{M} \bar{\bm{x}}_{i,N}^{[k]}\bm{D}_{i,j}^{[k]}-\bm{f}(\bar{\bm{x}}_{j,N}^{[k]},\bar{\bm{u}}_{j,N}^{[k]})|| \leq
\\ & || \sum_{i=0}^{M} \bar{\bm{x}}_{i,N}^{[k]}\bm{D}_{i,j}^{[k]} -\bm{\dot{x}}(t_{k,j}) || + 
\\ & \qquad ||\bm{f}(\bar{\bm{x}}_{j,N}^{[k]},\bar{\bm{u}}_{j,N}^{[k]}) - \bm{f}(\bm{x}(t_{k,j}),\bm{u}(t_{k,j}))|| \to 0  \, .
\end{split}
\end{equation*}
The convergence to zero follows from an application of the mean value theorem. The satisfaction of constraints Equations \eqref{eq:discequalityconstraint}, \eqref{eq:discinequalityconstraint} follow similarly. The satisfaction of Equations \eqref{eq:disccontinuityx} and \eqref{eq:disccontinuityt} follow from Equations \eqref{eq:tkj} and \eqref{eq:coeffdef}.

Before stating the main convergence result of this paper, the following is assumed. 

\begin{asm} \label{asm:convergence}
Let \(\{\bar{\bm{x}}_M^*,\bar{\bm{u}}_M^*,t_K^*\}\) be the sequence of optimal solutions to Problem \(P_M\), and let \(\bm{X}_M^*(t)\) and \(\bm{U}_M^*(t)\) be the control polygons defined by these solutions. There exist \(\bm{x}^{\infty}(t) \in \mathcal{C}^3, \, \bm{u}^{\infty}(t) \in \mathcal{C}^2\) on \([0,t_K^{\infty}]\) such that
\[
\lim_{K \to \infty, N \to \infty}(\bm{X}_M^*(t),\bm{U}_M^*(t),t_K^*) = (\bm{x}^{\infty}(t), \bm{u}^{\infty}(t),t_f^{\infty})
\]
for all \(t\in[0,t_f^\infty]\).
\end{asm}

\begin{thm} \label{thm:consistency}
Let \(\{\bm{x}_M^*(t),\bm{u}_M^*(t),t_K^*\}\) be the sequence of composite Bernstein polynomials obtained from the optimal solution to Problem \(P_M\), \(\bar{\bm{x}}_M^*\), \(\bar{\bm{u}}_M^*\) and final time \(t_K^*\), which satisfy Assumption \ref{asm:convergence}. 
Assume that Problem \ref{prob:continuous} has an optimal solution that satisfies \(\bm{x} \in \mathcal{C}^3\), \(\bm{u} \in \mathcal{C}^2\). Let the functions \(\bm{f}\), \(\bm{h}\) and \(F\) in Problem \ref{prob:continuous} be Lipschitz with respect to their arguments. 
Then, \(\{\bm{x}_M^*(t),\bm{u}_M^*(t),t_K^*\}\) converge to the optimal solution to Problem \(P\).
\end{thm}

\textbf{Proof:}
The proof is divided into three steps.

\textbf{Step $(1)$.} We prove that $(\bm{x}^\infty(t),\bm{u}^\infty(t),t_f^\infty)$ is a feasible solution to Problem $P$. We show by contradiction that $(\bm{x}^\infty(t),\bm{u}^\infty(t),t_f^\infty)$ satisfies the dynamic constraint of Problem $P$,
$
\bm{\dot{x}}^\infty(t)-\bm{f}(\bm{x}^\infty(t),\bm{u}^\infty(t)) = \bm{0} \, .
$
Assume that the above equality does not hold. Then, there exists $t'$, such that
\begin{equation} \label{eq:contradiction}
||\bm{\dot{x}}^\infty(t')-\bm{f}(\bm{x}^\infty(t'),\bm{u}^\infty(t'))|| > 0 \, .
\end{equation}
Since the nodes $\{t_{k,j}\}$ are dense in $[0,t_K]$, see Definition \ref{def:approximant}, there exist indexes \(k\) and \(j\) such that the following holds
$
\lim_{K,N \to \infty} t_{k,j} = t' .
$
Then, since $\bm{x}^\infty(t),\bm{u}^\infty(t) \in \mathcal{C}^2$, and from Assumption \ref{asm:convergence}, the left hand side of Equation \eqref{eq:contradiction} satisfies
\begin{equation*}
\begin{split}
& ||\bm{\dot{x}}^\infty(t')-\bm{f}(\bm{x}^\infty(t'),\bm{u}^\infty(t'))|| =
\\  & \qquad
\left\Vert \sum_{i=0}^{M} \bar{\bm{x}}_{i,N}^{[k]}\bm{D}_{i,j}^{[k]}-\bm{f}(\bar{\bm{x}}_{j,N}^{[k]},\bar{\bm{u}}_{j,N}^{[k]}) \right\Vert.
\end{split}
\end{equation*}
However, the dynamic constraint in Problem $P_M$ implies that
$
\lim_{K,N \to \infty} ||\sum_{i=0}^{M} \bar{\bm{x}}_{i,N}^{[k]}\bm{D}_{i,j}^{[k]}-\bm{f}(\bar{\bm{x}}_{j,N}^{[k]},\bar{\bm{u}}_{j,N}^{[k]})|| =  0,
$
which contradicts Equation \eqref{eq:contradiction}, proving that $\bm{x}^\infty(t)$ and $\bm{u}^\infty(t)$ satisfies the dynamic constraint in Equation~\eqref{eq:dynamicconstraint}. The equality and inequality constraints in Equations \eqref{eq:discequalityconstraint} and \eqref{eq:discinequalityconstraint} follow similarly.

\textbf{Step $(2)$.} We show that
$
\lim_{K,N \to \infty} I_M (\bx^*_M,\bu^*_M,\bar{\bm{t}}_K) = I(\bm{x}^\infty(t),\bm{u}^\infty(t),t_f)
$. I.e., we need to show that
\begin{equation} \label{eq:proofstep2a}
\begin{split}
& \lim_{K,N \to \infty} \sum_{k=1}^{K} w^{[k]}\sum_{j=0}^{N} F(\bm{x}_{j,N}^{[k]*},\bm{u}_{j,N}^{[k]*}) =  
\\ & \qquad  \int_0^{t_f^\infty} F(\bm{x}^\infty(t),\bm{u}^\infty(t))dt \, ,
\end{split}
\end{equation}
\begin{equation} \label{eq:proofstep2b}
\lim_{K,N \to \infty}  E(\bm{x}^*_M(0),\bm{x}^*_M(t_M)) = E(\bm{x}^\infty(0),\bm{x}^\infty(t_f)) \, .
\end{equation}
Using the Lipschitz assumption on $F$ and the continuity of $\bm{x}^\infty(t)$ and $\bm{u}^\infty(t)$, we get
\begin{equation*}
\begin{split} 
& \lim_{K,N \to \infty} \sum_{k=1}^{K} w^{[k]}\sum_{j=0}^{N} F(\bm{x}_{j,N}^{[k]*},\bm{u}_{j,N}^{[k]*}) = 
\\ & \qquad \int_0^{t_f} F(\bm{x}^\infty(t_{k,j}),\bm{u}^\infty(t_{k,j}))dt \, .
\end{split}
\end{equation*}
Finally, applying Lemma \ref{lem:integral}, the result in Equation \eqref{eq:proofstep2a} follows.
Similarly, using the Lipschitz assumption on $E$, one can show that Equation \eqref{eq:proofstep2b} holds.

\textbf{Step $(3)$.} We prove that $(\bm{x}^\infty(t),\bm{u}^\infty(t))$ is an optimal solution of Problem $P$, i.e.
$
I(\bm{x}^\infty(t),\bm{u}^\infty(t),t_f^\infty) = I(\bm{x}^*(t),\bm{u}^*(t),t_f) \, .
$
Let  $\tilde{\bx}_{M}$ and $\tilde{\bu}_{M}$ be the coefficients of the composite Bernstein polynomials approximating \(\bm{x}^*(t)\) and \(\bm{u}^*(t)\), respectively, and $\tilde{t_K}$ be the final time knot. 
An argument similar to the one in the proof of Step (2) yields
\begin{equation} \label{eq:convergenceoptimalcontinuous}
\lim_{K,N \to \infty} I_M(\tilde{\bx}_{M},\tilde{\bu}_{M},\tilde{t}_K) = I(\bm{x}^*(t),\bm{u}^*(t),t_f) \, .
\end{equation}
We note that
\begin{equation}
\begin{split}
& I(\bm{x}^*(t),\bm{u}^*(t),t_f) \leq I(\bm{x}^\infty(t),\bm{u}^\infty(t),t_f^{\infty})
\\ &  
= \lim_{K,N \to \infty} I_M(\bx^*_M,\bu^*_M,t_K) \leq \lim_{K,N \to \infty} I_M(\tilde{\bx}_{M},\tilde{\bu}_{M},\tilde{t}_K) \, ,
\end{split}
\end{equation}
which gives
$
I(\bm{x}^*(t),\bm{u}^*(t)) = I(\bm{x}^\infty(t),\bm{u}^\infty(t)) \, .
$

\begin{rem} The results presented above establish that solutions to the approximated Problem \ref{prob:discrete} exist and converge to the optimal solution of the original Problem \ref{prob:continuous}. 
The theorems are based on the convergence principles detailed in Lemmas \ref{lem:function}-\ref{lem:integral}, which form the basis for the assumptions in Theorems \ref{thm:feasibility} and \ref{thm:consistency}, where it is assumed that Problem \ref{prob:continuous} possesses solutions in \(\mathcal{C}^2\). These results could be extended by integrating insights from research on the composite Bernstein approximation of \(\mathcal{C}^0\) functions \cite{bojanic1989rate}, as well as functions with Lipschitz derivatives and H\"older continuous functions \cite{bojanic1989rate,mathe1999approximation,de1959degree}.
\end{rem}

\begin{rem} 
One distinct advantage of composite Bernstein polynomials lies in their capability and robustness in approximating discontinuous functions. This feature is critically relevant to optimal control, as many optimal control problems manifest bang-bang solutions. Traditional polynomial approximation methods, such as Fourier series, Lagrange interpolation, Hermite and Laguerre polynomial approximations, and Chebyshev or Legendre polynomials, often struggle with robustness when approximating step functions \cite{davis2022gibbs}. This typically manifests as oscillatory behavior near points of discontinuity, a phenomenon often referred to as the Gibbs phenomenon. However, Bernstein polynomials are notable for their immunity to this phenomenon \cite{gzyl2003approximation} and their convergence when approximating discontinuous functions with bounded variations \cite{bojanic1989rate}.
\end{rem}

\section{Knotting Method}\label{sec:KnottingMethod}
The collocation method presented in this paper can be used by simply predefining $K$ and optimizing the location of the knots. However, it might be more efficient to directly estimate the number of discontinuities and their locations to deduce an optimal $K$. When the number of discontinuities is known \emph{a priori}, it is apparent how many segments should be used to approximate the solution. However, in cases where the number of discontinuities is unknown, the number of segments required to approximate the solution is not clear. One such method used to identify these discontinuities is to evaluate the derivative of the control against some derivative threshold \cite{gong2008spectral}, chosen as a design variable by the user. Intuitively, this method identifies the location of sudden changes, i.e., the discontinuities, through the following steps: (1) Solve the NLP with $K=1$; (2) Calculate the derivative of the control input, i.e., $\dot{u}_N$; (3) Evaluate the elements of $\dot{u}_N$ against the threshold $D_{th}$; (4) For each element of $\dot{u}_N$ in violation of $D_{th}$, $K=K+1$.

Once the number of discontinuities is known, knots are placed at each detected discontinuity. Subsequently, the locations of the knots can be defined as a decision variable in the optimal control problem formulation, enabling the optimization algorithm to find their optimal placement.

\section{Numerical Results}\label{sec:NumericalResults}
In this section, two numerical examples are presented, highlighting the efficacy of using composite Bernstein polynomials to solve optimal control problems. These solutions were found using a constrained nonlinear optimization algorithm.
\subsection*{Bang-Bang Example \cite{tohidi2013efficient}}\label{subsec:example1}
Determine $y : [0,2] \to \mathbb{R}$ and $u:[0,2] \to \mathbb{R}$ that minimize
\begin{equation*}
    I(y(t),u(t)) = \int_0^2(3u(t)-2y(t))dt,
\end{equation*}
subject to
\begin{equation*}
    \dot{y}(t)=y(t)+u(t), \forall t \in [0,2],
\end{equation*}
\begin{equation*}
    y(0)=4, \, y(2)=39.392,
\end{equation*}
\begin{equation*}
    0 \le u(t) \le 2, \forall t \in [0,2].
\end{equation*}

\begin{figure*}
    \centering
    \begin{subfigure}{.32\textwidth}
        \centering
        \includegraphics[width=\textwidth,trim={.5cm 0cm .5cm .5cm},clip]{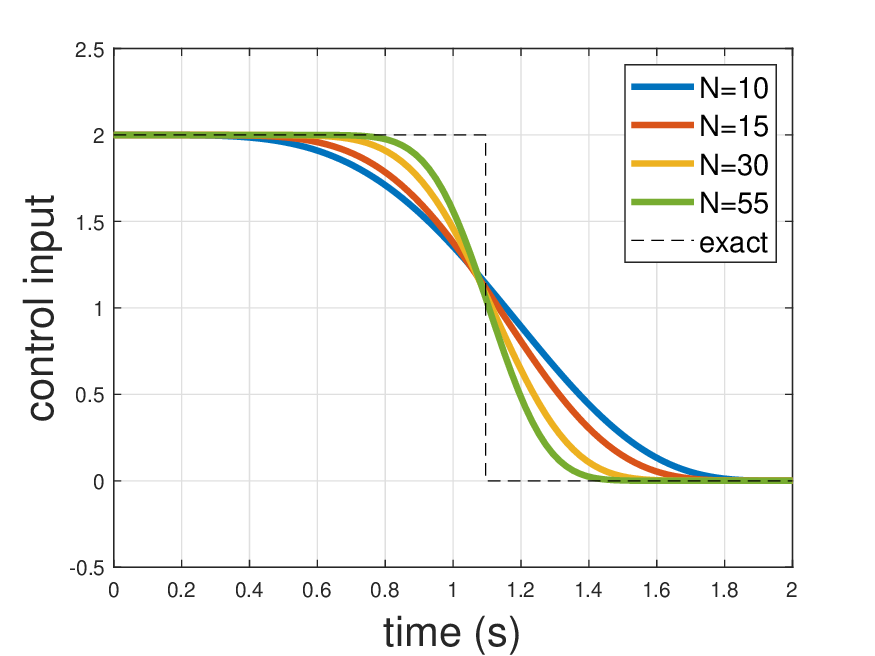}
        \label{fig:ex1_k1}
        \caption{}
    \end{subfigure}
    \begin{subfigure}{.32\textwidth}
        \centering
        \includegraphics[width=\textwidth,trim={.5cm 0cm .5cm .5cm},clip]{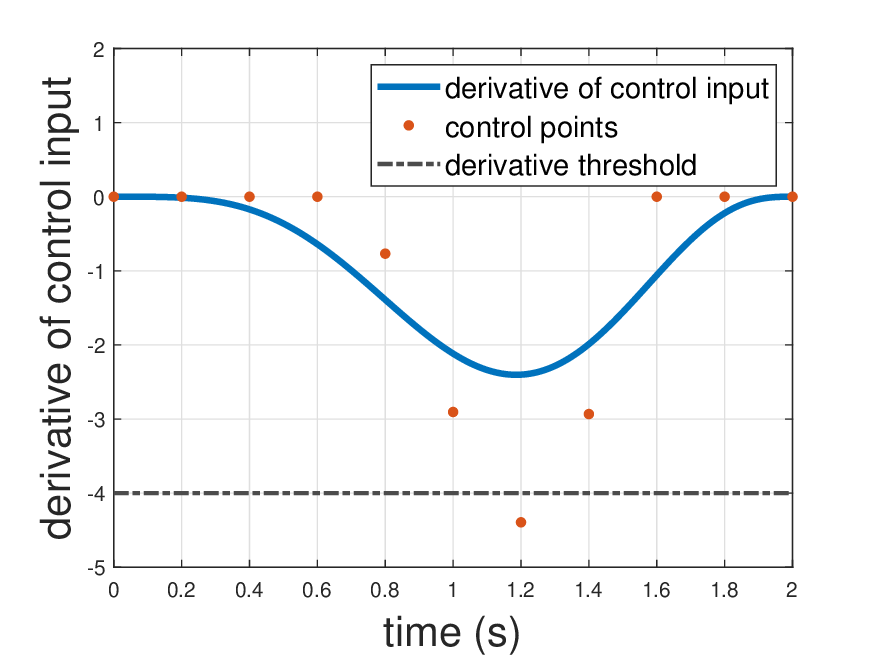}
        \label{fig:ex1_du}
        \caption{}
    \end{subfigure}
    \begin{subfigure}{.32\textwidth}
        \centering
        \includegraphics[width=\textwidth,trim={.5cm 0cm .5cm .5cm},clip]{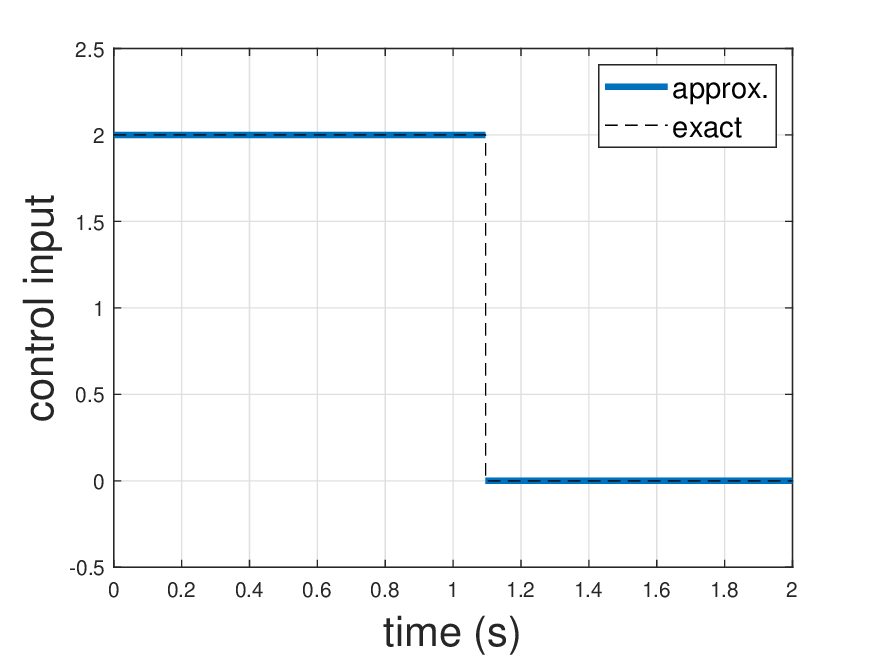}
        \label{fig:ex1_k2}
        \caption{}
    \end{subfigure}
    \caption{Solution to Example 1: (a) Approximation of a bang-bang input using a single Bernstein polynomial for orders $N=10,15,30,55$. (b) Derivative of the control input for the single Bernstein polynomial solution for $N=10$. (c) Approximation of a bang-bang input using a composite Bernstein polynomial consisting of two polynomials of order $N=10$.}
    \label{fig:ex1}
\end{figure*}

This problem is initially solved using the approximation method with one segment, i.e., the Bernstein approximation method. From Figure \ref{fig:ex1}(a), it is obvious that even though the approximation improves as $N$ increases, it remains inadequate for discontinuous solutions. The algorithm then evaluates the derivative of the control and compares it against some threshold to determine whether additional segments are required. A discontinuity is detected at $t=1.2s$, shown in Figure \ref{fig:ex1}(b), where the threshold is exceeded by one control point, indicating that two segments should be used to approximate the solution. The example is then solved with $K=2$, i.e., the composite Bernstein approximation method, with the results shown in Figure \ref{fig:ex1}(c). 

The composite Bernstein approximant, used with the presented knotting method, is able to detect the exact value of the discontinuity at $t=1.096s$, from the initial guess of $t=1.2s$. The evaluated cost functions of each scenario are shown in Table \ref{tab:ex1}. The cost of the analytical solution shown as a dotted line in Figures \ref{fig:ex1}(a) and \ref{fig:ex1}(c), is $J\approx-59.83$, the exact solution obtained with this method.

\begin{table}
    \centering
    \begin{tabular}{|c|c|c|}
    \hline
    Order & $K=1$ & $K=2$ \\
    \hline
    $N=10$ & $J\approx-59.35$ & $J\approx-59.83$ \\
    $N=15$ & $J\approx-59.50$ & - \\
    $N=30$ & $J\approx-59.67$ & - \\
    $N=55$ & $J\approx-59.74$ & - \\
    \hline
    \end{tabular}
    \caption{Bang-Bang Example Results: Evaluated cost function for each scenario.}
    \label{tab:ex1}
\end{table}



\subsection*{Multi-vehicle Motion Planning Problem} \label{sec:Application}
For this example, a multi-vehicle motion planning problem is formalized as an optimal control problem.
The scenario described below is that of an aircraft conducting a standard 45-degree traffic pattern entry which is subsequently alerted to a potential collision with an intruding aircraft, requiring replanning of the remaining trajectory to reach the runway at a free final time. 

Determine $x_1$, $x_2$, $x_3$, $x_4$, $x_5$, $u_1$, $u_2$, and $u_3$ that minimize
\begin{equation*}
    I(x(t),u(t))=t_f, 
\end{equation*}
subject to
\begin{equation*}
    \dot{x}_1(t)=u_1(t)\cos(x_4(t))\sin(x_5(t)), \forall t \in [0,t_f]
\end{equation*}
\begin{equation*}
    \dot{x}_2(t)=u_1(t)\sin(x_4(t))\cos(x_5(t)), \forall t \in [0,t_f]
\end{equation*}
\begin{equation*}
    \dot{x}_3(t)=u_1(t)\sin(x_5(t)), \forall t \in [0,t_f]
\end{equation*}
\begin{equation*}
    \dot{x}_4(t)=u_2(t), \, \dot{x}_5(t)=u_3(t), \forall t \in [0,t_f]
\end{equation*}
\begin{equation*}
    x_1(0)=0, \, x_2(0)=-4000,
\end{equation*}
\begin{equation*}
    x_3(0)=1000, \, x_4(0)= 180, \, x_5(0)=0,
\end{equation*}
\begin{equation*}
    x_1(t_f)=-3110.9, \, x_2(t_f)=0,
\end{equation*}
\begin{equation*}
    x_3(t_f)=259.2395, \, x_4(t_f)=0, \, x_5(t_f)=0,
\end{equation*}
\begin{equation*}
    0 \le x_3(t) \le 25000, \, -3 \le x_5(t) \le 3.5, \forall t\in[0,t_f]
\end{equation*}
\begin{equation*}
    100 \le u_1(t) \le 295, \, -3 \le u_2(t) \le 3, \forall t\in[0,t_f]
\end{equation*}
where $x_1$, $x_2$, and $x_3$ are the respective $[x,y,z]$ position of the aircraft in ft, $x_4$ is the heading in $\deg$, $x_5$ is the angle of attack in $\deg$, $u_1$ is the aircraft's velocity in $\frac{ft}{s}$, $u_2$ is the yaw rate in $\frac{\deg}{s}$, and $u_3$ is the pitch rate in $\frac{\deg}{s}$.
A separation variable of 500 ft is defined to specify the minimum safe distance between the aircraft and the intruding vehicle.

\begin{figure}[t]
    \centering
    \begin{subfigure}{.48\textwidth}
        \includegraphics[width=\textwidth,trim={0cm 2.5cm 0 3.25cm},clip]{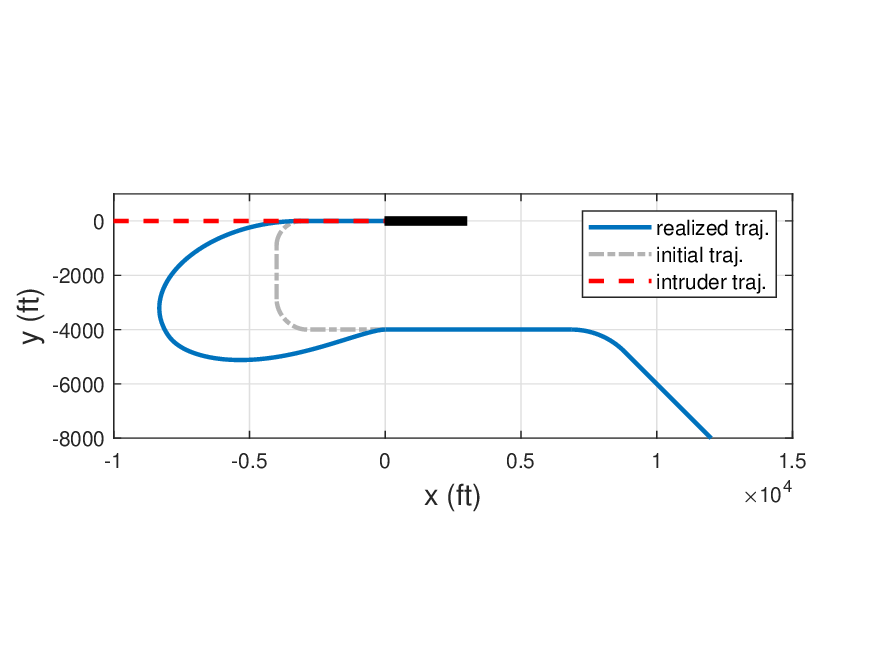}
        \label{fig:ex4_traj_2d}
        \caption{}
    \end{subfigure}
    \begin{subfigure}{.48\textwidth}
        \includegraphics[width=\textwidth,trim={0cm 0cm 0 0cm},clip]{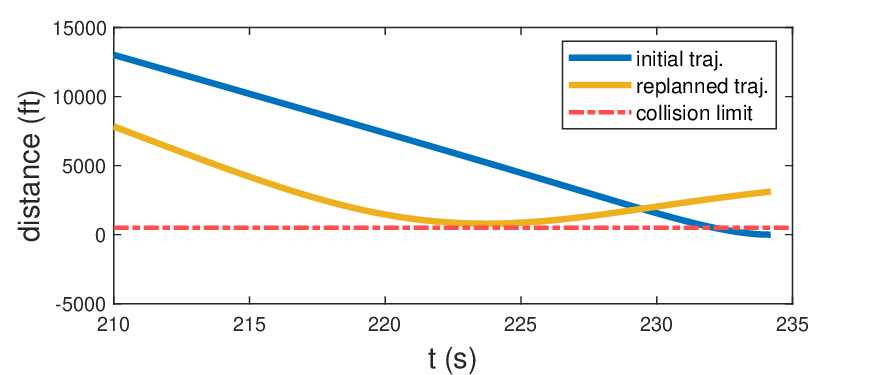}
        \label{fig:ex4_coll}
        \caption{}
    \end{subfigure}
    \caption{(a) Aerial map of runway with aircraft following 45$\degree$ entry until an intruding aircraft is detected and the initial trajectory is replanned to avoid collision. (b) Distance between the initial and replanned trajectory of the aircraft, and the trajectory of the intruder for the motion planning example, with a separation variable of 500 ft.}
    \label{fig:ex4}
\end{figure}

The initial, replanned, and intruder trajectories are shown in Figure \ref{fig:ex4}(a). A potential collision is detected when the vehicle is at roughly $[0,-4000,1000]$ ft, thus replanning occurs. With the replanned trajectory, the vehicle evades the intruder, avoiding collision. This maneuver is shown in Figure \ref{fig:ex4}(b). Flight trajectories, like the initial trajectory shown in Figure \ref{fig:ex4}(a), exhibit explicit geometric patterns contingent upon vehicle position and velocity. This example highlights the ability of composite Bernstein polynomials to precisely capture specific flight patterns, which a single Bernstein polynomial would struggle to represent at low orders.


\section{Conclusion}\label{sec:Conclusion}
In this paper, we introduced composite Bernstein polynomials as a means to solve optimal control problems as a direct method, via direct approximation of the continuous optimal control problem into a discrete solution. This method builds off of the favorable properties of Bernstein polynomials for motion planning, and extends this work to composite Bernstein polynomials, allowing for faster convergence and more accurate approximations of discontinuous solutions. Additionally, rigorous analyses of the convergence properties of composite Bernstein polynomials were presented. Two numerical examples and their solutions via this method were exhibited, including a motion planning problem, demonstrating the efficacy of using composite Bernstein polynomials to solve complex optimal control problems.

\appendices
\section{Proof of Lemma \ref{lem:function}} 
\label{app.lem:function}
First, notice that 
\begin{equation}\label{eq:tau}
\Vert {x}_M(t) - {x}(t) \Vert \leq \max_{ k = 1, \ldots, K } \max_{\tau \in [t_k,t_{k+1}]} \Vert {x}_N^{[k]}(\tau) - {x}(\tau) \Vert
\end{equation}
for all \(k \in  1 , \ldots , K \) and for all \(\tau \in [t_{k-1}, t_{k}] \), the following equality holds:
\begin{equation*}
\begin{split}  
\frac{(\tau-t_{k-1})(t_{k}-\tau)}{N} \dot{b}_{j,N}^{[k]}(\tau) =  \\
\left( t_{k-1} + j \frac{t_{k}-t_{k-1}}{N} - \tau \right) b_{j,N}^{[k]}(\tau) \, .	 
\end{split} 
\end{equation*}
Let $[t_0 , \ldots , t_k ; x]$ denote the $k$-th order divided difference of $x(\tau)$ at the points $t_0 , \ldots , t_k$. Then, the \emph{generalized} Stancu's remainder formula \cite{floater2005convergence} is derived as follows:
\begin{equation*}
\begin{split}  
& x_N^{[k]}(\tau) - x(\tau) \\
& = \sum_{j=0}^N b_{j,N}^{[k]}(\tau) \left(x\left( t_{k-1} + j \frac{t_k-t_{k-1}}{N} \right) - x(\tau) \right)  \\
&  = \sum_{j=0}^N b_{j,N}^{[k]}(\tau) \left[t_{k-1} + j \frac{t_k-t_{k-1}}{N},\tau;x\right] \\  
\end{split}
\end{equation*}
\begin{equation*}
\begin{split}
& \qquad \left( t_{k-1} +   j \frac{t_k-t_{k-1}}{N} - \tau \right) \\
& = \frac{(\tau-t_{k-1})(t_k-\tau)}{N} \sum_{j=0}^N  \dot{b}_{j,N}^{[k]} (\tau) \\
& \left[t_{k-1} + j \frac{t_k-t_{k-1}}{N},\tau;x\right] \\
& = \frac{(\tau-t_{k-1})(t_k-\tau)}{(t_{k}-t_{k-1})} \sum_{j=0}^N  \left[t_{k-1} + j \frac{t_k-t_{k-1}}{N},\tau;x\right] \\
& \left( b_{j-1,N-1}^{[k]}(\tau) -b_{j,N-1}^{[k]}(\tau) \right)
\end{split} 
\end{equation*}
\begin{equation} \label{app.eq:stancugeneralized}
\begin{split}
& \quad = \frac{(\tau-t_{k-1})(t_k-\tau)}{(t_k-t_{k-1})} \sum_{j=0}^{N-1}  \left( \left[t_{k-1} +  j\frac{t_k-t_{k-1}}{N},\tau;x\right] 
\right. \\ & \qquad \left. 
- \left[t_{k-1} +  (j+1)\frac{t_k-t_{k-1}}{N},\tau;x\right] \right) b_{j,N-1}^{[k]}(\tau) \\
& \quad
= \frac{(\tau-t_{k-1})(t_k-\tau)}{N} \sum_{j=0}^{N-1} \left[t_{k-1} +  j\frac{t_k-t_{k-1}}{N}, t_{k-1} 
\right. \\ & \qquad \left.
+ (j+1)\frac{t_k-t_{k-1}}{N},\tau;x\right] b_{j,N-1}^{[k]}(\tau) \,
\end{split} 
\end{equation}
noticing that $b_{j,N-1}^{[k]}(\tau)=0$ for $0 < j < N-1$, Equation \eqref{app.eq:stancugeneralized} simplifies to
\begin{equation*}
\begin{split}
    |x_N^{[k]}(\tau) - x(\tau)|\le \frac{(\tau-t_{k-1})(t_k-\tau)}{N}\frac{1}{2}  \max_{\tau \in [t_{k-1} \, , \, t_k]} |\ddot x (\tau)| \\
    \le \frac{(\frac{t_k-t_{k-1}}{2})^2}{N}\frac{1}{2}\max_{\tau \in [t_{k-1} \, , \, t_k]}|\ddot x (\tau)|.
\end{split}
\end{equation*}
Thus, the following result holds, and proves Lemma \ref{lem:function}:
$$
| x_N^{[k]}(\tau) - x(\tau) | \leq \frac{C_t^2}{8K^2N} \max_{\tau \in [t_{k-1} \, , \, t_k]} |\ddot x (\tau)| \, .
$$


\section{Proof of Lemma \ref{lem:derivative}} 
\label{app.lem:derivative}
Let us define the following operator:
\begin{equation}\label{eq:Boperator}
\begin{split}
    B_{n,s,m}x(\tau)=\sum^{n-s}_{j=0}\left[t_{k-1}+j\frac{t_k-t_{k-1}}{N},...,t_{k-1}+ \right. \\ \left. (j+s)\frac{t_k-t_{k-1}}{N},\underbrace{\tau,...,\tau}_m;x\right]b_{j,n-s}
\end{split}
\end{equation}
for all $k\in {1,...,K}$ and for all $\tau\in [t_{k-1},t_k]$.
Then, Equation \eqref{app.eq:stancugeneralized} can be rewritten as follows:
\begin{equation}\label{eq:preleibniz}
    x^{[k]}_N(\tau)-x(\tau)=\frac{(\tau-t_{k-1})(t_k-\tau)}{N}B_{N,1,1}x(\tau).
\end{equation}
Differentiation of Equation \eqref{eq:preleibniz} using the Leibniz rule gives
\begin{equation}\label{eq:derivativeError}
\begin{split}
    & x^{[k](r)}_N(\tau)-x^{(r)}(\tau)\\ & = \sum^r_{q=0}\binom{r}{q}\frac{d^q}{d\tau^q}\left(\frac{(\tau-t_{k-1})(t_k-\tau)}{N}\right) B_{N,1,1}^{(r-q)}x(\tau) \\ & =
    \frac{(\tau-t_{k-1})(t_k-\tau)}{N}B_{N,1,1}^{(r)}x(\tau)\\ & +\frac{r(t_k-2\tau+t_{k-1})}{N}B_{N,1,1}^{(r-1)}x(\tau) -\frac{r(r-1)}{N}B_{N,1,1}^{(r-2)}x(\tau).
\end{split}
\end{equation}
Now we investigate the derivatives of $B_{N,1,1}x(\tau)$. By using the following relationship
\begin{equation*}
\begin{split}
    & \frac{d^r}{d\tau^r}\left[t_{k-1}+j\frac{t_k-t_{k-1}}{N},t_{k-1}+(j+1)\frac{t_k-t_{k-1}}{N},\tau;x\right] \\
    & = r!\left[t_{k-1}+j\frac{t_k-t_{k-1}}{N},t_{k-1}
    \right. \\ & \quad \left. +(j+1)\frac{t_k-t_{k-1}}{N},\underbrace{\tau,...,\tau}_{r+1};x\right].
\end{split}
\end{equation*}
Differentiation of Equation \eqref{eq:Boperator} with $s=m=1$ gives
\begin{equation} \label{eq:BoperatorDerivative}
\begin{split}
    (B_{N,1,1}x)^{(r)}(\tau)
    = \sum^{N-1}_{j=0}\sum^r_{q=0}\binom{r}{q}(r-q)!\left.(\left.[t_{k-1}+\right.\right. \\
    j\frac{t_k-t_{k-1}}{N}
    ,t_{k-1}\left.\left.
    +(j+1)\frac{t_k-t_{k-1}}{N},\underbrace{\tau,...,\tau}_{r-q+1};x\right.]\right.)b_{j,N-1}^{(q)}(\tau) \\
    =r!\sum^r_{q=0}\frac{(N-1)...(N-q)}{q!(t_k-t_{k-1})^q}\sum^{N-q-1}_{j=0}  (\Delta^q[t_{k-1} +j\frac{t_k-t_{k-1}}{N},\\t_{k-1}+(j+1)\frac{t_k-t_{k-1}}{N},\underbrace{\tau,...,\tau}_{r-q+1};x]) \times b_{j,N-q-1}(\tau)
\end{split}
\end{equation}
where $\Delta$ is the forward difference operator w.r.t. $j$. Notice that
\begin{equation*}
\begin{split}
    \Delta[t_{k-1}+j\frac{t_k-t_{k-1}}{N},t_{k-1} 
    +(j+1)\frac{t_k-t_{k-1}}{N}, \\
    \tau,...,\tau;x]
    =[t_{k-1}+(j+1)\frac{t_k-t_{k-1}}{N},t_{k-1}\\
    +(j+2)\frac{t_k-t_{k-1}}{N},\tau,...,\tau;x]-[t_{k-1} +(j+1) \\
    \frac{t_k-t_{k-1}}{N},t_{k-1}+(j+2)\frac{t_k-t_{k-1}}{N},\tau,...,\tau;x] \\
    =\frac{2(t_k-t_{k-1})}{N}[t_{k-1}+j\frac{t_k-t_{k-1}}{N},t_{k-1} \\
    +(j+1)\frac{t_k-t_{k-1}}{N},t_{k-1}+(j+2)\frac{t_k-t_{k-1}}{N},\tau,...,\tau;x]
\end{split}
\end{equation*}
and continuing to apply $\Delta$ implies
\begin{equation*}
\begin{split}
    &\Delta^q\left[t_{k-1}+j\frac{t_k-t_{k-1}}{N},t_{k-1} \right. \\ 
    & \quad \left. +(j+1)\frac{t_k-t_{k-1}}{N},\tau,...,\tau;x\right]
    \\ & =\frac{(q+1)!(t_k-t_{k-1})^q}{N^q}[t_{k-1}+j\frac{t_k-t_{k-1}}{N},...,t_{k-1} 
    \\&\quad +(j+q+1)\frac{t_k-t_{k-1}}{N},\tau,...,\tau;x].
\end{split}
\end{equation*}
Substituting the last result into Equation \eqref{eq:BoperatorDerivative} and replacing $q$ with $q-1$ gives
\begin{equation*}
\begin{split}
    (B_{N,1,1}x)^{(r)}(\tau)=r!\sum^{r+1}_{q=1}q\frac{(N-1)...(N-q+1)}{N^{q-1}} \\
    \times\sum^{N-q-1}_{j=0} ([t_{k-1}+j\frac{t_k-t_{k-1}}{N},...,t_{k-1} \\
    +(j+q)\frac{t_k-t_{k-1}}{N},\underbrace{\tau,...,\tau}_{r-q+2};x])b_{j,N-q}(\tau) 
\end{split}
\end{equation*}
\begin{equation*}
\begin{split}
    =r!\sum^{r+1}_{q=1}q\frac{(N-1)...(N-q+1)}{N^{q-1}}B_{N,q,r-q+2}(\tau).
\end{split}
\end{equation*}
From \cite{floater2005convergence} and the previous equation, we can conclude the following
\begin{equation}
    \Vert(B_{N,1,1}x)^{(r)}(\tau)\Vert\le r!\sum^{r+1}_{q=1}q\frac{\Vert x^{(r+2)}\Vert}{(r+2)!}\le\frac{\Vert x^{(r+2)}\Vert}{2}.
\end{equation}
Recalling Equation \eqref{eq:derivativeError}, we get
\begin{equation}\label{eq:derivativeError2}
\begin{split}
    & |x_N^{[k](r)}(\tau)-x^{(r)}(\tau)| \\ &\le \frac{1}{2N}((\tau-t_{k-1})(t_k-\tau)\Vert x^{(r+2)}(\tau)\Vert \\
    & \quad +r(t_k-2\tau+t_{k-1})\Vert x^{(r+1)}(\tau)\Vert+r(r-1)\Vert x^{(r)}(\tau)\Vert).
\end{split}
\end{equation}
Noticing that
\begin{equation*}
\begin{split}
    \Vert x^{(r+1)}(\tau)\Vert = \left\Vert\int^{\tau}_{t_{k-1}} x^{(r+2)}(t)dt + x^{(r+1)}(t_{k-1})\right\Vert \\
    \le \int_{t_{k-1}}^{t_k}\Vert x^{(r+2)}(\tau)\Vert d\tau + \Vert x^{(r+1)}(t_{k-1})\Vert \\
    \le (t_k-t_{k-1})\Vert x^{(r+2)}(\tau)\Vert + \Vert x^{(r+1)}(t_{k-1})\Vert \\
    \le \frac{C_t}{K}\Vert x^{(r+2)}(\tau)\Vert + \Vert x^{(r+1)}(t_{k-1})\Vert,
\end{split}
\end{equation*}
and that
\begin{equation*}
\begin{split}
    \Vert x^{(r)}(\tau)\Vert \le \frac{C_t}{K}\Vert x^{(r+1)}(\tau)\Vert + \Vert x^{(r)}(t_{k-1})\Vert \\
    \le \frac{C_t^2}{K^2}\Vert x^{(r+2)}(\tau)\Vert+\frac{C_t}{K}\Vert x^{(r+1)}(t_{k-1})\Vert+\Vert x^{(r)}(t_{k-1})\Vert,
\end{split}
\end{equation*}
Equation \eqref{eq:derivativeError2} expands to
\begin{equation*}
\begin{split}
    |x_N^{[k](r)}(\tau)-x^{(r)}(\tau)| \le \frac{C_t^2}{2K^2N}\Vert x^{(r+2)}(\tau) \\
    +\frac{rC_t}{2KN}\left(\frac{C_t}{K}\Vert x^{(r+2)}(\tau)+\Vert x^{(r+1)}(t_{k-1})\Vert\right) \\
    +\frac{r(r-1)}{2N}\left(\frac{C_t^2}{K^2}\Vert x^{(r+2)}(\tau)\Vert+\frac{C_t}{K}\Vert x^{(r+1)}(t_{k-1})\Vert \right. \\ \left. +\Vert x^r(t_{k-1})\Vert\right.\bigg).
\end{split}
\end{equation*}
Then, Lemma \ref{lem:derivative} follows with Equation \eqref{eq:lemma2results}.


\section{Proof of Lemma \ref{lem:integral}}
\label{app.lem:integral}
Using Lemma \ref{lem:function}, and the triangle inequality for integrals, we note that
\begin{equation*}
\begin{split}
    \left\Vert\int^{t_k}_{t_{k-1}}x_N^{[k]}(\tau)d\tau-\int^{t_k}_{t_{k-1}}x(\tau)d\tau\right\Vert \\
    \le\int^{t_k}_{t_{k-1}}\Vert x_N^{[k]}(\tau)-x(\tau)\Vert d\tau \\
    \le\int^{t_k}_{t_{k-1}}\frac{C_t^2}{8K^2N}\max_{\tau\in[t_{k-1},t_k]}|\ddot{x}(\tau)|d\tau. \\
\end{split}
\end{equation*}
Thus
\begin{equation*}
\begin{split}
    \left\Vert\int^{t_f}_{t_o}x_M(t)dt-\int^{t_f}_{t_o}x(t)dt\right\Vert \\
    \le\sum_{k=1}^K\int^{t_k}_{t_{k-1}}\frac{C_t^2}{8K^2N}\max_{\tau\in[t_{k-1},t_k],k=1,...,K}|\ddot{x}(\tau)|d\tau \\
    = \frac{C_t^3}{8K^2N}\max_{\tau\in[t_{k-1},t_k],k=1,...,K}|\ddot{x}(\tau)|
\end{split}
\end{equation*}
which proves Lemma \ref{lem:integral}.

\bibliographystyle{IEEEtran}
\bibliography{refs}

\begin{thebibliography}{10}
\providecommand{\url}[1]{#1}
\csname url@samestyle\endcsname
\providecommand{\newblock}{\relax}
\providecommand{\bibinfo}[2]{#2}
\providecommand{\BIBentrySTDinterwordspacing}{\spaceskip=0pt\relax}
\providecommand{\BIBentryALTinterwordstretchfactor}{4}
\providecommand{\BIBentryALTinterwordspacing}{\spaceskip=\fontdimen2\font plus
\BIBentryALTinterwordstretchfactor\fontdimen3\font minus \fontdimen4\font\relax}
\providecommand{\BIBforeignlanguage}[2]{{%
\expandafter\ifx\csname l@#1\endcsname\relax
\typeout{** WARNING: IEEEtran.bst: No hyphenation pattern has been}%
\typeout{** loaded for the language `#1'. Using the pattern for}%
\typeout{** the default language instead.}%
\else
\language=\csname l@#1\endcsname
\fi
#2}}
\providecommand{\BIBdecl}{\relax}
\BIBdecl

\bibitem{karpenko2012first}
M.~Karpenko, S.~Bhatt, N.~Bedrossian, A.~Fleming, and I.~Ross, ``First flight results on time-optimal spacecraft slews,'' \emph{Journal of Guidance, Control, and Dynamics}, vol.~35, no.~2, pp. 367--376, 2012.

\bibitem{ross2006pseudospectral}
I.~Ross, P.~Sekhavat, A.~Fleming, Q.~Gong, and W.~Kang, ``Pseudospectral feedback control: foundations, examples and experimental results,'' in \emph{AIAA guidance, navigation, and control conference and exhibit}, 2006.

\bibitem{elnagar1995pseudospectral}
G.~Elnagar, M.~A. Kazemi, and M.~Razzaghi, ``The pseudospectral {L}egendre method for discretizing optimal control problems,'' \emph{IEEE transactions on Automatic Control}, vol.~40, no.~10, pp. 1793--1796, 1995.

\bibitem{ross2004legendre}
I.~M. Ross and F.~Fahroo, ``Legendre pseudospectral approximations of optimal control problems,'' in \emph{New trends in nonlinear dynamics and control and their applications}.\hskip 1em plus 0.5em minus 0.4em\relax Springer, 2004, pp. 327--342.

\bibitem{benson2005gauss}
D.~Benson, ``A gauss pseudospectral transcription for optimal control,'' Ph.D. dissertation, Massachusetts Institute of Technology, 2005.

\bibitem{huntington2007advancement}
G.~T. Huntington, ``Advancement and analysis of a gauss pseudospectral transcription for optimal control problems,'' Ph.D. dissertation, Massachusetts Institute of Technology, Department of Aeronautics and Astronautics, 2007.

\bibitem{kameswaran2008convergence}
S.~Kameswaran and L.~T. Biegler, ``Convergence rates for direct transcription of optimal control problems using collocation at {R}adau points,'' \emph{Computational Optimization and Applications}, vol.~41, no.~1, pp. 81--126, 2008.

\bibitem{garg2011direct}
D.~Garg, M.~A. Patterson, C.~Francolin, C.~L. Darby, G.~T. Huntington, W.~W. Hager, and A.~V. Rao, ``Direct trajectory optimization and costate estimation of finite-horizon and infinite-horizon optimal control problems using a {R}adau pseudospectral method,'' \emph{Computational Optimization and Applications}, vol.~49, pp. 335--358, 2011.

\bibitem{ross2012review}
I.~M. Ross and M.~Karpenko, ``A review of pseudospectral optimal control: from theory to flight,'' \emph{Annual Reviews in Control}, vol.~36, no.~2, pp. 182--197, 2012.

\bibitem{sarra2009edge}
S.~A. Sarra, ``Edge detection free postprocessing for pseudospectral approximations,'' \emph{Journal of Scientific Computing}, vol.~41, no.~1, pp. 49--61, 2009.

\bibitem{kielas2022bernstein}
C.~Kielas-Jensen, V.~Cichella, T.~Berry, I.~Kaminer, C.~Walton, and A.~Pascoal, ``Bernstein polynomial-based method for solving optimal trajectory generation problems,'' \emph{Sensors}, vol.~22, no.~5, p. 1869, 2022.

\bibitem{cichella2018bernstein}
V.~Cichella, I.~Kaminer, C.~Walton, N.~Hovakimyan, and A.~Pascoal, ``Bernstein approximation of optimal control problems,'' \emph{arXiv preprint arXiv:1812.06132}, 2018.

\bibitem{bojanic1989rate}
R.~Bojanic and F.~Cheng, ``Rate of convergence of {B}ernstein polynomials for functions with derivatives of bounded variation,'' \emph{Journal of Mathematical Analysis and Applications}, vol. 141, no.~1, pp. 136--151, 1989.

\bibitem{mathe1999approximation}
P.~Math{\'e}, ``Approximation of {H}{\"o}lder continuous functions by {B}ernstein polynomials,'' \emph{The American mathematical monthly}, vol. 106, no.~6, pp. 568--574, 1999.

\bibitem{de1959degree}
K.~de~Leeuw, ``On the degree of approximation by {B}ernstein polynomials,'' \emph{Journal d’Analyse Math{\'e}matique}, vol.~7, no.~1, pp. 89--104, 1959.

\bibitem{davis2022gibbs}
J.~M. Davis and P.~Hagelstein, ``Gibbs phenomena for some classical orthogonal polynomials,'' \emph{Journal of Mathematical Analysis and Applications}, vol. 505, no.~1, p. 125574, 2022.

\bibitem{gzyl2003approximation}
H.~Gzyl and J.~L. Palacios, ``On the approximation properties of {B}ernstein polynomials via probabilistic tools,'' \emph{Bolet{\i}n de la Asociaci{\'o}n Matem{\'a}tica Venezolana}, vol.~10, no.~1, pp. 5--13, 2003.

\bibitem{gong2008spectral}
Q.~Gong, F.~Fahroo, and I.~M. Ross, ``Spectral algorithm for pseudospectral methods in optimal control,'' \emph{Journal of Guidance, Control, and Dynamics}, vol.~31, no.~3, pp. 460--471, 2008.

\bibitem{tohidi2013efficient}
E.~Tohidi, A.~Pasban, A.~Kilicman, and S.~L. Noghabi, ``An efficient pseudospectral method for solving a class of nonlinear optimal control problems,'' in \emph{Abstract and Applied Analysis}, vol. 2013.\hskip 1em plus 0.5em minus 0.4em\relax Hindawi, 2013.

\bibitem{floater2005convergence}
M.~S. Floater, ``On the convergence of derivatives of {B}ernstein approximation,'' \emph{Journal of Approximation Theory}, vol. 134, no.~1, pp. 130--135, 2005.

\end{thebibliography}

\end{document}